# PROPORTIONAL FAIRNESS AND ITS RELATIONSHIP WITH MULTI-CLASS QUEUEING NETWORKS

By N. S. Walton[1]

*University of Cambridge*

We consider multi-class single-server queueing networks that have a product form stationary distribution. A new limit result proves a sequence of such networks converges weakly to a stochastic flow level model. The stochastic flow level model found is insensitive. A large deviation principle for the stationary distribution of these multi-class queueing networks is also found. Its rate function has a dual form that coincides with proportional fairness. We then give the first rigorous proof that the stationary throughput of a multi-class single-server queueing network converges to a proportionally fair allocation.

This work combines classical queueing networks with more recent work on stochastic flow level models and proportional fairness. One could view these seemingly different models as the same system described at different levels of granularity: a microscopic, queueing level description; a macroscopic, flow level description and a teleological, optimization description.

**1. Introduction.** In this paper we form descriptions of multi-class single-server queueing networks at different levels of granularity. Similar descriptions of electrical networks have been well studied and provide a good analogue of the results proven in this paper.

One could form a Markov chain model of electrons in an electrical network. At this first level, one explicitly describes the location of particles. Transitions within the network occur rapidly, so perhaps it is more natural to consider results like Ohm's law and Kirchhoff's law which are concerned with the current flowing through the network. At this second level, one considers the average flow of particles through the network. An electrical network also minimizes energy dissipation, as described by Thomson's principle

Received September 2008; revised April 2009.
[1]Supported in part by EPSRC Grant GR/586266/01.
*AMS 2000 subject classifications.* Primary 60K25, 60K30; secondary 90K15, 68K20.
*Key words and phrases.* Multi-class queueing network, proportional fairness, bandwidth sharing, stochastic flow level model, insensitivity, product form stationary distribution, proportionally fair, state space collapse.

(primal) and Dirichlet's principle (dual). At this third and final level, one considers the network to be acting as an optimizer. For further discussion see Kelly [14], Section 2, and Doyle and Snell [7], Section 1.

Just as we consider electrons, current and energy minimization in an electrical network, respectively, in a model of a packet switched network we consider packets, bandwidth and utility optimization. Following Shah and Wischik [24], we will use the terms microscopic, macroscopic and teleological to refer to these different descriptions of our network. *Microscopic* refers to a detailed description of a network's state though gives little insight into overall dependence. *Macroscopic* refers to an averaged view of the original network. *Teleological* considers a further abstracted view where the network can be seen to be acting globally as an optimizer.

In this paper our microscopic model will be product form single-server multi-class queueing networks [11]; our macroscopic model will be a specific stochastic flow level model [2, 3, 8, 18, 19, 21] and our teleological model will be the proportionally fair optimization problem [15, 25]. The main results of this paper are concerned with forming rigorous connections between these different models.

The first result of this paper is concerned with connecting our microscopic model to our macroscopic model. Stochastic flow level models for an intuitive model of document transfer across a packet switched network. Despite this no rigorous convergence proof has been constructed to justify a stochastic flow level model as the limit of a packet switching network. In Theorem 3.1 we construct a proof to address this issue. We view a sequence of multi-class queueing networks as a simplistic model of document transfer across a packet switching network, and we prove weak convergence in the Skorohod topology of these networks to a specific stochastic flow level model. We call the resulting stochastic flow level model the spinning network. The spinning network was first considered by Massoulié [21], Section 3.4, and formed the first insensitive stochastic flow level model. Insensitivity results on this model are given by Proutière [21] and Bonald and Proutière [3].

The second result of this paper is concerned with connecting both our microscopic and macroscopic models to our teleological model. In particular, in Theorems 6.1 and 7.2, we give the first rigorous proofs of a mathematical relationship between multi-class networks of single-server queues and proportional fairness. An argument justifying such a relationship has been made by Schweitzer [23], Kelly [13] and Massoulié and Roberts [19]. This noteworthy argument, presented in Section 5, considers the constraints that a queueing network may impose on transfer rates. Here we take a different approach and provide a rigorous proof using large deviations and convex duality. The use of large deviations is motivated by the relationship between balanced fairness and proportional fairness found by Massoulié [17]. In addition, Pittel [20] has considered the large deviations of multi-class queueing networks but does not derive proportional fairness.

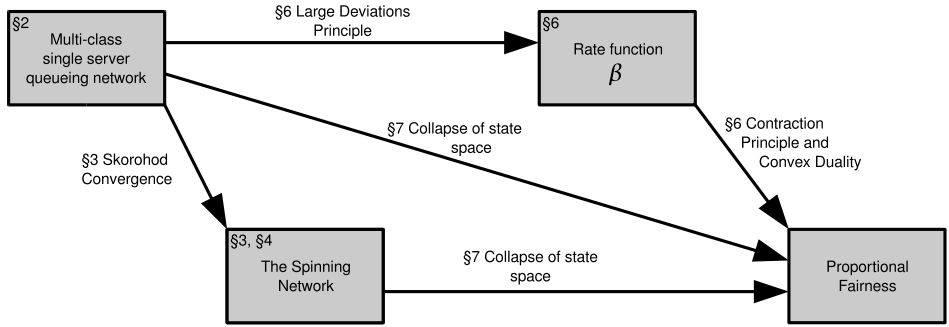

Fig. 1. *Structure of the paper's results.*

To prove Theorem 6.1 we consider a large deviations principle for the stationary distribution of a multi-class single-server queueing network and find a rate function $\beta(\cdot)$. Applying the Contraction Principle, we gain a new rate function $\alpha(\cdot)$, expressed as a convex optimization problem. In primal form, $\alpha(\cdot)$ is interpreted as minimizing entropy subject to constraints. We find that the dual form of $\alpha(\cdot)$ is, up to a constant, the proportionally fair optimization problem. These arguments give Theorem 6.1. With this we are able to prove Theorem 7.2. Theorem 7.2 states that the stationary throughput of a closed multi-class queueing network converges to a proportionally fair allocation as the number of customers (or packets) is increased in proportion to some fixed vector. Proportionally fair optimization occurs because the states of the queueing network collapse to a set of entropy minimizing states with proportionally fair throughput. These results are analogous to the heavy-traffic notion of state space collapse [4, 16, 24]. These results emphasize a large deviations duality between network state and network flow. Please see Figure 1 for a diagram of the structure of these results and the sections that they are contained in.

In addition, the multi-class queueing networks considered in this paper are known to be quasi-reversible and thus have a product form equilibrium distribution [11]. The equivalence between reversibility and insensitivity is well studied [3, 21, 22, 26]. Our macroscopic model, the spinning network, inherits the reversibility property from these multi-class queueing networks and thus, as has been observed by Massoulié, the spinning network is insensitive. Thus for multi-class queueing networks we find at each different level of granularity a different feature of the system can be observed: product form at the microscopic level; insensitivity at the macroscopic level and proportional fairness at the teleological level.

Observations on insensitivity and the product form stationary behavior of proportional fairness have previously been made. By considering balanced fairness, the connection between insensitivity and proportional fairness has been given by Massoulié [17]. Also, motivated by diffusion approximation behavior in queueing networks, Kang, Kelly, Lee and Williams [9, 10, 16] considered heavy-traffic approximations of stochastic flow level models operating under proportional fairness and found them to have a product form stationary distribution.

A final observation is that the multi-class queueing networks considered here have no prescribed optimization structure. Thus it is surprising to see that asymptotically these networks are implicitly solving a utility optimization problem.

1.1. *Organization.* The sections of the paper are structured as follows. In the next subsection we introduce notation, define proportionally fair bandwidth allocations and the proportionally fair optimization problem. In Section 2 we define the multi-class queueing networks considered in this paper and state results on their stationary distribution. In Section 3 we introduce stochastic flow level models, define the spinning network and state Theorem 3.1. In Section 4 we discuss the reversibility, insensitivity, stability and stationary behavior of the spinning network. In Section 5 we present the argument given by Schweitzer [23], Kelly [13] and Massoulié and Roberts [19] relating proportional fairness and multi-class queueing networks. In Section 6 we consider the large deviations behavior of the multi-class single-server queueing networks presented in Section 2 and prove Theorem 6.1. In Section 7 we discuss the collapse in state space brought about by large deviations principle Theorem 6.1 and prove convergence of stationary throughput in Theorem 7.2. In the Appendix we give the proof of Theorem 3.1 and prove additional lemmas and propositions from within the text.

1.2. *Notation: Network structure and proportional fairness.* We let a finite set $\mathcal{J}$ index the set of *queues* in a network. Let $J = |\mathcal{J}|$. A *route* through the network is a nonempty set of queues. Let $\mathcal{I} \subset 2^{\mathcal{J}}$ be the set of routes. Let $I = |\mathcal{I}|$. For each route $i = \{j_1^i, \ldots, j_{k_i}^i\} \in \mathcal{I}$, we associate an order $(j_1^i, \ldots, j_{k_i}^i)$. Also we define the set of queue-route incidences, $\mathcal{K} := \{(j,i) : i \in \mathcal{I}, j \in \mathcal{J}, j \in i\}$ and let $K = |\mathcal{K}|$. We will view our multi-class queueing network model as transferring a number of documents across the different routes of the network. In Sections 2, 3, 4 and Appendix A.1, the vector $n = (n_i : i \in \mathcal{I}) \in \mathbb{Z}_+^I$ will be used to refer to the number of documents in transfer across the routes of the network. When referring to large deviations characteristics, in Sections 6 and 7, $n = (n_i : i \in \mathcal{I}) \in \mathbb{R}_+^I$ will be used to refer to the proportion of documents in transfer across routes. We will consider the documents in our multi-class queueing network to be transferred by packets which will traverse the network. Each document on each route will have only one packet in transfer across the network at any point in time. Thus, given our description of $n$, the total number (or proportion) of packets in transfer across route $i$ will be $n_i$. We also consider the vector $m = (m_{ji} : (j,i) \in \mathcal{K}) \in \mathbb{R}_+^K$. In Sections 2, 3, 4 and Appendix A.1, $m_{ji} \in \mathbb{Z}_+$ will be used to refer to the number of packets in transfer across route $i$ that are at queue $j$. Similarly in Sections 6 and 7, $m_{ji} \in \mathbb{R}_+$ will be used to refer to the proportion of the packets in transfer, that are on route $i$ and in queue $j$. Because each packet in transfer will correspond to a document we have that

$$n_i = \sum_{j \,:\, j \in i} m_{ji} \qquad \forall i \in \mathcal{I}.$$

We define the number (or proportion) of packets at a queue to be

$$m_j := \sum_{i \,:\, j \in i} m_{ji} \qquad \forall j \in \mathcal{J}.$$

For each $n \in \mathbb{Z}_+^I$ we define $S(n) = \{m \in \mathbb{Z}_+^K : \sum_{j \,:\, j \in i} m_{ji} = n_i \ \forall i \in \mathcal{I}\}$, the set of queue states achievable given the number of documents in transfer.

With each queue $j \in \mathcal{J}$ we associate a service capacity $C_j$. A *bandwidth allocation* is a vector $\Lambda(n) = (\Lambda_i(n) : i \in \mathcal{I}) \in \mathbb{R}_+^I$ for each $n \in \mathbb{Z}_+^I$. A bandwidth allocation is said to be *feasible*, if $\forall n \in \mathbb{Z}_+^I$

$$\sum_{i \,:\, j \in i} \Lambda_i(n) \leq C_j \qquad \forall j \in \mathcal{J}.$$

A bandwidth allocation corresponds to the rate documents are transferred across their route given the number of documents in transfer. Bandwidth allocations form an abstraction of the stationary transfer rate achieved by rate control algorithms used, for example, in the internet.

A bandwidth allocation $\Lambda^{\text{PF}}(n)$ is *proportionally fair* [15] if $\forall n \in \mathbb{R}_+^I$, $\Lambda_i^{\text{PF}}(n) = 0$ when $n_i = 0$ and $\Lambda^{\text{PF}}(n)$ solves

$$\text{(1.1)} \qquad \text{maximize} \sum_{i \in \mathcal{I}} n_i \log \Lambda_i,$$

$$\text{(1.2)} \qquad \text{subject to} \sum_{i \,:\, j \in i} \Lambda_i \leq C_j \qquad \forall j \in \mathcal{J},$$

$$\text{(1.3)} \qquad \text{over } \Lambda_i \geq 0 \qquad \forall i \in \mathcal{I}.$$

For all vectors $x \in \mathbb{R}^d$ we define $\lfloor x \rfloor := (\lfloor x_1 \rfloor, \ldots, \lfloor x_d \rfloor)$, the lower integer part of each component. Unless stated otherwise $\|x\| = \max_{d'} |x_{d'}|$, the supremum norm. Finally, we will define that for each $m \in \mathbb{Z}_+^K$

$$\binom{m_j}{m_{ji} : i \ni j} = \frac{m_j!}{\prod_{i \,:\, j \in i} (m_{ji}!)}.$$

**2. A microscopic queueing model.** In this section we introduce our microscopic model. The queueing networks considered here are exactly the networks with fixed service capacity described in Sections 3.1 and 3.4 of Kelly [11]. We interpret these multi-class queueing networks as a Markov chain model of document transfer across a packet switching network. This interpretation has previously been considered by Massoulié and Roberts [19] and Bonald and Proutière [3]. Documents wishing to be transferred across a network are broken into a number of packets. The document's packets are then sent across the network one by one, so that a new packet is sent into the network once the packet in the network has completed its route. We define our queueing model in the next three paragraphs and call it an *open multi-class queueing network with spinning*.

We consider a network of queues indexed by the set $\mathcal{J}$ queues process packets. Each packet moves along a fixed route from the set $\mathcal{I}$. Each queue $j \in \mathcal{J}$ may store an infinite number of packets and has a fixed service capacity $C_j \in (0, \infty)$. Each packet at a queue has a position in that queue, for example, if there are $m_j$ packets at queue $j$ then these packets are stored in positions $k = 1, \ldots, m_j$. The total service capacity of a queue is then divided between the different packets at the queue. Each queue operates under a service discipline that cannot discriminate between the routes used by its packets. More explicitly there exists a function $\gamma_j(k, m_j)$ that gives the proportion of service devoted to the packet in position $k$ in queue $j$ when there are $m_j$ packets at queue $j$. As $\gamma_j$ represents a proportion

$$(2.1) \qquad \sum_{k=1}^{m_j} \gamma_j(k, m_j) = 1 \qquad \forall m_j > 0.$$

Similarly when joining a queue a packet may only chose its position as a function of the number of packets at that queue. These service disciplines are described in Kelly [11], pages 58–60, and, for example, include first in first out, last come first served, processor sharing, and symmetric queues.

Documents for transfer on route $i \in \mathcal{I}$ arrive as a Poisson process of rate $\nu_i$. Each route $i$ document consists of a discrete number of packets. This number is independent, finite mean and with distribution equal to random variable $X_i$. We let $\mu_i = (\mathbb{E}X_i)^{-1}$ and let $\rho_i = \frac{\nu_i}{\mu_i}$ for all $i \in \mathcal{I}$. Each packet has an independent, exponentially distributed mean 1 service requirement at each queue.

Consider route $i \in \mathcal{I}$ with route order $(j_1^i, \ldots, j_{k_i}^i)$. If we wish to transfer a document across route $i$, a packet is sent along route $i$. It will first join queue $j_1^i$. For $k = 1, \ldots, k_i - 1$ on departing queue $j_k^i$ a packet will join queue $j_{k+1}^i$. When a packet leaves its final queue $j(i) = j_{k_i}^i$, a new packet is sent along route $i$ until all packets in the document are transferred.

Equivalently, one could think of each document being transferred by a single packet which repeats its route with some probability. This probability only depends on the number of repetitions of its route the packet has made so far. This interpretation motivates the use of the word "spinning." These models are equivalent to the networks with fixed service capacity considered in Kelly [11], Section 3.1. All the results above and all the proportional fairness results in Sections 6 and 7 apply to this case.

We could explicitly describe the state of this network by recording the position of each packet at each queue, the route used by these packets, the number of repetitions such packets have made on their route and the total number of repetitions these packets must make. As noted in [11], the stochastic process recording this information is a Markov chain. We will not be interested in this explicit description. We will be interested in simpler quantities, namely the number of documents in transfer across routes, the number of packets in transfer on each route at each queue and the throughput of packets of each route at each queue. In general, the processes associated with these quantities will not be Markov.

We now consider the stationary distribution of this model. The following result is a direct consequence of the Theorem 3.1 of Kelly [11].

PROPOSITION 2.1. *An open multi-class queueing network with spinning is ergodic if and only if*

$$(2.2) \qquad \sum_{i\,:\,j\in i} \rho_i < C_j \qquad \forall j \in \mathcal{J}.$$

*When ergodic, $M = (M_{ji}\colon (j,i) \in \mathcal{K})$, the number of packets in transfer across each route at each queue, has stationary distribution,*

$$(2.3) \qquad \mathbb{P}(M=m) = B^{-1} \prod_{j \in \mathcal{J}} \left( \binom{m_j}{m_{ji}\,:\,i \ni j} \prod_{i\,:\,j \in i} \left(\frac{\rho_i}{C_j}\right)^{m_{ji}} \right)$$

*for each $m \in \mathbb{Z}_+^K$, where*

$$(2.4) \qquad B := \prod_{j \in \mathcal{J}} \left( \frac{C_j}{C_j - \sum_{i \ni j} \rho_i} \right).$$

PROOF. Allow the state of a packet to be given by the packet's route, the packet's position in its current queue, the total number of repetitions the packet makes of its route and the number of repetitions currently made. From this we have a Markov chain of the form described in Section 3.1 of [11]. Applying Theorem 3.1 of [11] to find the stationary distribution and summing over the correct states gains the result. □

The next two corollaries are an immediate consequence of this result.

COROLLARY 2.1. *$N = (N_i : i \in \mathcal{I})$ the number of documents in transfer has stationary distribution*

$$\mathbb{P}(N = n) = \frac{B_n}{B} \prod_{i \in \mathcal{I}} \rho_i^{n_i} \qquad \forall n \in \mathbb{Z}_+^I, \tag{2.5}$$

*where we define*

$$B_n := \sum_{m \in S(n)} \prod_{j \in \mathcal{J}} \left( \binom{m_j}{m_{ji} : i \ni j} \prod_{i \,:\, j \in i} \left( \frac{1}{C_j} \right)^{m_{ji}} \right) \qquad \forall n \in \mathbb{Z}_+^I. \tag{2.6}$$

A *closed multi-class queueing network* behaves as an open multi-class queueing network except that document arrivals and departures are forbidden (see [11]). So, the network behaves as if there is a fixed number of infinitely large documents in transfer. We assume throughout this paper:

ASSUMPTION 1. Consider the Markov chain description of a closed queueing network that records the position of each packet at each queue and each packet's route [11], Section 3.4. Given the number of packets on each route, consider the set of all possible states of the Markov chain. We assume this set of states is irreducible.

This assumption excludes reducibility issues which can only occur in closed queueing networks where a queue serves a single deterministically chosen packet. It is worth noting that if Assumption 1 is broken then there need not be a unique stationary distribution or a unique stationary throughput for the closed queueing network.

COROLLARY 2.2. *For a closed multi-class queueing network with $n \in \mathbb{Z}_+^I$ documents in transfer across routes, the number of packets in transfer of each route at each queue has stationary distribution*

$$\mathbb{P}_n(M = m) = B_n^{-1} \prod_{j \in \mathcal{J}} \left( \binom{m_j}{m_{ji} : i \ni j} \prod_{i \,:\, j \in i} \left( \frac{1}{C_j} \right)^{m_{ji}} \right) \tag{2.7}$$

*for each $m \in \mathbb{Z}_+^K$, where $B_n$ is defined by (2.6).*

Finally we can characterize the stationary throughput of these closed multi-class queueing networks.

COROLLARY 2.3. *For a closed multi-class queueing network with $n \in \mathbb{Z}_+^I$ documents in transfer across routes and with $n_i > 0$, the stationary throughput of route $i$ packets at queue $j \in i$ is*

$$\Lambda_i^{\text{SN}}(n) := \frac{B_{n-e_i}}{B_n},$$

*where $B_n$ is defined by (2.6) and $e_i$ is the $i$th unit vector in $\mathbb{R}_+^I$.*

PROOF. The probability the network is in state $m \in \mathbb{Z}_+^K$ is given by (2.7). Given the network is in state $m$, by Corollary 3.4 of [11], the probability in queue $j$ the packet position $k \in \{1, \ldots, m_j\}$ is traversing route $i$ is $\frac{m_{ji}}{m_j}$. The throughput of the packet in position $k$ of queue $j$ is $\gamma_j(k, m_j) C_j$. Thus the stationary throughput of the network is

$$\sum_{\substack{m \in S(n): \\ m_j > 0}} \sum_{k=1}^{m_j} \gamma_j(k, m_j) C_j \frac{m_{ji}}{m_j} \frac{1}{B_n} \prod_{l \in \mathcal{J}} \left( \binom{m_l}{m_{lr}: r \ni l} \prod_{r: l \in r} \left(\frac{1}{C_l}\right)^{m_{lr}} \right)$$

$$= \sum_{\substack{m \in S(n): \\ m_j > 0}} C_j \frac{m_{ji}}{m_j} \frac{1}{B_n} \prod_{l \in \mathcal{J}} \left( \binom{m_l}{m_{lr}: r \ni l} \prod_{r: l \in r} \left(\frac{1}{C_l}\right)^{m_{lr}} \right)$$

$$= \sum_{m' \in S(n-e_i)} \frac{1}{B_n} \prod_{l \in \mathcal{J}} \left( \binom{m'_l}{m'_{lr}: r \ni l} \prod_{r: l \in r} \left(\frac{1}{C_l}\right)^{m'_{lr}} \right) = \frac{B_{n-e_i}}{B_n}.$$

We used (2.1) in the first inequality; in the second we cancelled terms and substituted $m'_{lr} = m_{lr} - 1$ if $(l, r) = (j, i)$ and $m'_{lr} = m_{lr}$ otherwise. □

We will define $\Lambda_i^{\text{SN}}(n)$ as a bandwidth allocation in the next section.

**3. Limit to a macroscopic model.** For an electrical network at the macroscopic level, we considered the current through the network and not the explicit behavior of individual electrons. Similarly for a packet switched network, we may wish to consider the rate documents are transferred through the network and not the explicit location of packets. Stochastic flow level models provide such a model of document transfer as they do not explicitly consider packets. In this section we justify how a series of multi-class queueing networks with spinning converges in the Skorohod topology to a stochastic flow level model. Thus we limit from a model where documents are transferred by sending discrete packets to a model where documents are transferred at a dynamic, elastic rate.

A key quantity in this analysis will be the stationary throughput of a closed multi-class queueing network.

DEFINITION 1 (Spinning allocation). For all $n \in \mathbb{Z}_+^I$, the spinning allocation, denoted $\Lambda^{\text{SN}}(n) = (\Lambda_i^{\text{SN}}(n): i \in \mathcal{I})$, is the stationary throughput of packets on each route of a closed multi-class queueing network with $n$ documents in transfer. More explicitly from Corollary 2.3 we know that, $\forall i \in \mathcal{I}$ and $n \in \mathbb{Z}_+^I$

$$(3.1) \quad \Lambda_i^{\text{SN}}(n) = \begin{cases} \dfrac{B_{n-e_i}}{B_n}, & \text{if } n_i > 0, \\ 0, & \text{otherwise,} \end{cases}$$

where $B_n$ is defined by (2.6).

Under the name "The store-forward allocation," the spinning allocation is cited by Proutière [21] as the first insensitive bandwidth allocation. The definition of this bandwidth allocation from multi-class queueing networks is due to Massoulié [21], Section 3.4. Our macroscopic model of interest will be the following.

DEFINITION 2 (Stochastic flow level model). A stochastic flow level model operating under bandwidth allocation $\Lambda(\cdot)$ is a continuous-time Markov chain on $\mathbb{Z}_+^I$ with rates

$$(3.2) \quad q(n,n') = \begin{cases} \nu_i, & \text{if } n' = n + e_i, \\ \mu_i \Lambda_i(n), & \text{if } n' = n - e_i \text{ and } n_i > 0, \\ 0, & \text{otherwise}, \end{cases}$$

$e_i$ is the $i$th unit vector in $\mathbb{Z}_+^I$.

DEFINITION 3 (Spinning network). The spinning network is the stochastic flow level model operating under the spinning allocation.

Stochastic flow level models were first considered by Roberts and Massoulié [18]. This model can be interpreted as follows. Documents wishing to be transferred across route $i$ arrive as a Poisson process of rate $\nu_i$. These documents are assumed to have a size that is independent and exponentially distributed with mean $\mu_i^{-1}$. If currently the number of documents in transfer across routes is given by vector $n \in \mathbb{Z}_+^I$, then each document on route $i$ is transferred at rate $\frac{\Lambda_i(n)}{n_i}$. Documents are then processed at this rate until there is a change in the network's state, either by a document transfer being completed and thus leaving the network, or by a document arrival occurring. Thanks to the memoryless property of our process we need not record residual document sizes when an arrival or departure event occurs.

We now introduce the sequence of multi-class queueing networks which we will limit to form our macroscopic model. Consider a sequence of open multi-class queueing networks with spinning, $\{(M^{(c)}(t) : t \in \mathbb{R}_+)\}_{c \in \mathbb{N}}$. These networks have the same routing structure and document arrival processes as described in the last section. In this section, we assume for simplicity that each queue is processor sharing. Thus $M^{(c)}$ is a Markov chain. We increase the rate at which packets are transferred through the network. In the $c$th network each queue $j$ operates at service rate $cC_j$. We also increase the size of documents, so as not to increase the rate that documents are transferred through the network. We assume route $i$ documents are geometrically distributed with parameter $\mu_i/c$. We also let $N^{(c)}$ be the stochastic process for the number of documents of each route in transfer, that is, $\forall i \in \mathcal{I}$, $\forall t \in \mathbb{R}_+$

$$N_i^{(c)}(t) = \sum_{j\,:\,j \in i} M_{ji}^{(c)}(t).$$

Let us consider intuitively how these networks limit as $c \to \infty$. Note transitions of packets between queues occur at times of order $O(\frac{1}{c})$. Thus the number of packets sent by a route $i$ document in transfer per unit time is $O(c)$. The probability that a packet sent is the final packet is $\frac{\mu_i}{c}$, so the time until this document is transferred is $\frac{\mu_i}{c} O(c) = O(1)$. Thus there is a separation of time scales between document transfer and packet transfer. Between arrival and departure times the network behaves as a closed queueing network. By the ergodic theorem, for large $c$, this closed queueing network will behave close to its stationary distribution. Thus between arrival and departure times documents experience a transfer rate determined by the stationary throughput of a closed queueing network. Noting Definition 1 and Corollary 2.3, we have defined this transfer rate to be the spinning allocation. Also, the increased rate of packet transfer and the geometric number of packets in a document suggests an exponential distribution limit for the size of documents. Thus it seems plausible that a stochastic flow level model would be the limit of these queueing networks and that $\forall n \in \mathbb{Z}_+^I$, the rate of transfer, would be determined by the spinning allocation.

We will prove this assertion is correct. The formal statement of this convergence result is the following theorem. For our multi-class queueing network the theorem rigorises the separation of times scales assumption in Massoulié and Roberts [18, 19]. Here we define the Skorohod topology with any norm on $\mathbb{R}^I$.

THEOREM 3.1. *For each $c \in \mathbb{N}$, take an open multi-class queueing network with spinning $M^{(c)}$ as described above. We assume queues are processor sharing. Let $N^{(\infty)}$ denote the number of documents in transfer in the spinning network, (3.1) and (3.2). If*

$$N^{(c)}(0) \Rightarrow N^{(\infty)}(0) \qquad as\ c \to \infty,$$

*then, in the Skorohod topology on interval $[0, 1]$,*

$$N^{(c)} \Rightarrow N^{(\infty)} \qquad as\ c \to \infty.$$

Due to its technical nature this result is proven in the Appendix. The proof uses a coupling argument. The key idea in this proof is to let the internal behavior of each queueing network between arrival and departure times be governed by the same process while allowing the number of packets sent before a departure to converge almost surely.

REMARK 1. One could use a similar model of a multi-class queueing network with spinning where the number of packets sent is not geometrically distributed and perform this limit. In this way one would model the transfer of documents of any positive distribution. Over $c \in \mathbb{N}$ and these different document size distributions, the stationary distribution of $N^{(c)}$ will be unchanged provided the mean document size is scaled so that $\mu^{(c)} = (c\mu_i : i \in \mathcal{I})$.

REMARK 2. The processor sharing assumption is not needed in Theorem 3.1. In general only Assumption 1 will be needed. We can prove Theorem 3.1 when queues are not processor sharing by replacing $M^{(c)}$ with the explicit Markov chain description of the queueing network outlined in Section 2.

**4. A macroscopic stochastic flow level model.** In the last section we introduced stochastic flow level models and justified how one such model, the spinning network, formed a macroscopic model of the queueing networks considered in Section 2. Now we will discuss some properties of the spinning network. In particular, we show its stationary distribution to be insensitive to different document size distributions. Though we first note:

LEMMA 4.1. *The spinning allocation is a feasible bandwidth allocation.*

PROOF. For processor sharing queues, $\Lambda_i^{\mathrm{SN}}(n) = \mathbb{E}_n[\frac{M_{ji}}{M_j} C_j \mathbb{I}[M_j > 0]]$, so

$$\sum_{i:j \in i} \Lambda_i^{\mathrm{SN}}(n) = \sum_{i:j \in i} \mathbb{E}_n \left[ \frac{M_{ji}}{M_j} C_j \mathbb{I}[M_j > 0] \right] = C_j \mathbb{P}_n(M_j > 0) \leq C_j. \qquad \square$$

We can extend the definition of a stochastic flow level model so that the sizes of incoming documents are independent and of any positive distribution. Information on residual document sizes would be needed for such processes to be Markov. Given this extension, a stochastic flow level model with mean document sizes given by $(\mu_i^{-1} : i \in \mathcal{I})$ is *insensitive* if the stationary distribution for the number of documents in transfer is the same as all other stochastic flow level models with the same mean document sizes.

The stationary distribution of an open multi-class queueing network with spinning (2.3) depends on the distribution of the number of packets in a document only through mean document size $(\mu_i^{-1} : i \in \mathcal{I})$. In this sense an open multi-class queueing network with spinning is insensitive. By the same scaling in Theorem 3.1 we could increase the network's service capacity and limit the discrete document size distribution to approximate continuous document size distributions. Under this scaling the stationary distribution (2.3) still depends on the distribution of the number of packets in a document only through parameters $(\mu_i : i \in \mathcal{I})$. Thus given Theorem 3.1, it is reasonable to think that its limit, the spinning network, would be insensitive.

Bonald and Proutière [3, 21] found that key results on the insensitivity of stochastic flow level models are a consequence of existing results on the insensitivity and reversibility of Whittle networks [22, 26].

PROPOSITION 4.1 (Bonald and Proutière [3]). *An ergodic stochastic flow level model operating under bandwidth allocation $\Lambda(\cdot)$ is insensitive if and only if it is reversible, that is, there exists function $\Phi : \mathbb{Z}^I \to \mathbb{R}_+$ with $\Phi(0) = 1$, $\Phi(n) = 0 \ \forall n \notin \mathbb{Z}_+^I$ and*

$$\Lambda_i(n) = \frac{\Phi(n - e_i)}{\Phi(n)} \qquad \forall n \in \mathbb{Z}_+^I, i \in \mathcal{I},$$

*moreover,*

(4.1) $$\pi(n) = \Phi(n) \prod_{i \in \mathcal{I}} \rho_i^{n_i},$$

*forms an invariant measure for the number of documents in transfer.*

We now find the stationary distribution of the spinning network and show it to be insensitive; this fact has been observed by Proutière [21] and Bonald and Proutière [3].

PROPOSITION 4.2. *The spinning network is reversible and insensitive to document size distributions. The spinning network is ergodic if and only if*

(4.2) $$\sum_{i \,:\, j \in i} \rho_i < C_j \qquad \forall j \in \mathcal{J}.$$

*The spinning network has the same stationary distribution as the number of documents in transfer in an open multi-class queueing network with spinning, that is distribution, (2.5).*

PROOF. Insensitivity and the reversible property are an immediate consequence of Proposition 4.1 and Definition 1. By Proposition 4.1

$$\pi(n) := B_n \prod_{i \in \mathcal{I}} \rho_i^{n_i} \qquad \forall n \in \mathbb{Z}_+^I,$$

is an invariant measure. The sum of $\pi(\cdot)$ over all states is finite if and only if the stability condition (4.2) holds. When finite this sum equals $B$ given by (2.4) thus giving stationary distribution (2.5). □

**5. Relating multi-class queueing networks to proportional fairness.** In 1979, Schweitzer [23] studied approximations of closed multi-class queueing networks and considered how asymptotic conditions on such networks might satisfy the Kuhn–Tucker conditions for proportionally fair optimization. In 1989, Kelly [13] studied approximations of closed queueing networks and by an analogous analysis considered a similar optimization formulation. In 1999, Massoulié and Roberts [19] studied a fluid-type queueing model and used these same Kuhn–Tucker conditions to deduce proportional fairness. To develop intuition and to motivate Sections 6 and 7, we present the argument used in these three papers.

As given in Section 2, consider a closed multi-class queueing network with $n_i$ documents in transfer on each route $i \in \mathcal{I}$. Let $q_j$ be the mean sojourn time of a packet at queue $j$; let $\bar{m}_{ji}$ be the mean number of route $i$ packets in transfer at queue $j$, and let $\Lambda_i$ be the mean sending rate of route $i$ packets into the network. By Little's law,

$$\Lambda_i q_j = \bar{m}_{ji} \qquad \forall (j,i) \in \mathcal{K}. \tag{5.1}$$

Summing over $j \in i$ and rearranging gives

$$\frac{n_i}{\Lambda_i} - \sum_{j \in i} q_j = 0 \qquad \forall i \in \mathcal{I}. \tag{5.2}$$

Since queues are stable we know

$$\sum_{i \,:\, j \in i} \Lambda_i \leq C_j \qquad \forall j \in \mathcal{J}. \tag{5.3}$$

One can imagine if equality (5.3) is strict then $q_j \approx 0$. Thus approximately

$$q_j \left( C_j - \sum_{i \,:\, j \in i} \Lambda_i \right) = 0 \qquad \forall j \in \mathcal{J}. \tag{5.4}$$

Also

$$q_j \geq 0 \quad \forall j \in \mathcal{J} \quad \text{and} \quad \Lambda_i \geq 0 \quad \forall i \in \mathcal{I}. \tag{5.5}$$

Interpreting $(q_j : j \in \mathcal{J})$ as Lagrange multipliers, (5.2)–(5.5) are precisely the Kuhn–Tucker conditions for the proportionally fair optimization problem

$$\max_{\Lambda \in \mathbb{R}_+^I} \sum_{i \in \mathcal{I}} n_i \log \Lambda_i \quad \text{subject to} \quad \sum_{i \,:\, j \in i} \Lambda_i \leq C_j \qquad \forall j \in \mathcal{J}.$$

So from this one can deduce that $\Lambda_i = \Lambda_i^{\mathrm{PF}}(n)$ $\forall i \in \mathcal{I}$.

To make this argument we assumed that the sojourn times of packets did not depend on the route used and that complementary slackness condition (5.4) held. Neither of these conditions need be true in general. In fact, from Corollary 2.3, know that $\Lambda_i = \Lambda_i^{\mathrm{SN}}(n)$, the spinning allocation. In general $\Lambda_i^{\mathrm{SN}}(n) \neq \Lambda_i^{\mathrm{PF}}(n)$. Even so it is reasonable to assume $\Lambda_i^{\mathrm{SN}}(n) \approx \Lambda_i^{\mathrm{PF}}(n)$.

In the following two sections we rigorously prove a relationship between multi-class queueing networks and proportional fairness. We consider a multi-class network of single-server queues, as described in Section 2. We let the number of documents in transfer get large but in proportion to some fixed vector $n \in \mathbb{R}_+^I$. We show that these multi-class queueing networks asymptotically allocate service across routes as a proportionally fair optimizer. For example, we will prove that for all $i \in \mathcal{I}$

$$\Lambda_i^{\mathrm{SN}}(\lfloor hn \rfloor) \xrightarrow[h \to \infty]{} \Lambda_i^{\mathrm{PF}}(n).$$

**6. Limit to a teleological description.** In the Introduction we noted how minimizing energy dissipation gave an optimization description of an electrical network. In this section we wish to justify how proportional fairness provides an optimization description for the open multi-class queueing networks discussed in Section 2.

To do this we allow the number of documents in transfer to be large and in proportion to some fixed vector. The main result in this section is Theorem 6.1 where we prove a large deviation principle for stochastic models with stationary distribution (2.5). This stationary distribution includes the number of documents in transfer for all open multi-class queueing networks discussed in Section 2 and the spinning network with any finite mean document size distribution. This large deviations approach for queueing networks is similar to that given by Pittel [20], although Pittel does not consider a relationship with proportional fairness in his analysis. In addition, as we consider the large deviations of an insensitive stochastic flow level model this approach is also similar to that taken by Massoulié [17] for balanced fairness. The large deviation rate function found in Theorem 6.1 is

$$(6.1) \quad \alpha_\rho(n) = \max_{\Lambda \in \mathbb{R}_+^I} \sum_{i\,:\,n_i>0} n_i \log \frac{\Lambda_i}{\rho_i} \quad \text{subject to} \quad \sum_{i\,:\,j \in i} \Lambda_i \leq C_j \quad \forall j \in \mathcal{J}.$$

When optimizing the above expression we can express the $\rho_i$ terms as additive constants. Thus the argument maximizing this optimization problem is the proportionally fair allocation $\Lambda^{\mathrm{PF}}(n)$. From this we see that these queueing models are related to proportionally fair optimization.

To prove Theorem 6.1, first we prove a large deviation principle for the stationary distribution (2.3). Stirling's formula finds a rate function $\beta_\rho(\cdot)$. Applying the contraction mapping principle gives the large deviation principle for the number of documents in transfer and finds $\alpha_\rho(\cdot)$ expressed as the primal of a convex optimization problem. We calculate the dual of this optimization problem and find it to be of the form of (6.1).

We start by finding rate function $\beta_\rho(\cdot)$.

LEMMA 6.1. *Suppose $M$ is a random variable in $\mathbb{Z}_+^K$ with distribution (2.3). If we take a vector $m \in \mathbb{R}_+^K$ and take $\{d^{(h)}\}_{h \in \mathbb{N}}$ a sequence of vectors in $\mathbb{R}^K$ such that $hm + d^{(h)} \in \mathbb{Z}_+^K$ and $\sup_h |d^{(h)}| < \infty$ then*

$$\lim_{h \to \infty} \frac{1}{h} \log \mathbb{P}(M = hm + d^{(h)}) = -\beta_\rho(m),$$

*where we define*

(6.2) $$\beta_\rho(m) := \sum_{\substack{(j,i) \in \mathcal{K}: \\ m_j > 0}} m_{ji} \log \frac{m_{ji} C_j}{m_j \rho_i}.$$

PROOF. For all $j \in \mathcal{J}$, define $d_j^{(h)} = \sum_{i : j \in i} d_{ji}^{(h)}$. By Stirling's formula,

$$\lim_{h \to \infty} \frac{1}{h} \log \mathbb{P}(M = hm + d^{(h)})$$

$$= \lim_{h \to \infty} \frac{1}{h} \Bigg[ \sum_{j \in \mathcal{J}} \log(hm_j + d_j^{(h)})! - \sum_{(j,i) \in \mathcal{K}} \log(hm_{ji} + d_{ji}^{(h)})!$$

$$+ \sum_{(j,i) \in \mathcal{K}} (hm_{ji} + d_{ji}^{(h)}) \log \frac{\rho_i}{C_j} \Bigg]$$

$$= \lim_{h \to \infty} \frac{1}{h} \Bigg[ \sum_{\substack{j \in \mathcal{J}: \\ m_j > 0}} ((hm_j + d_j^{(h)}) \log(hm_j + d_j^{(h)}) - (hm_j + d_j^{(h)}))$$

$$- \sum_{\substack{(j,i) \in \mathcal{K}: \\ m_{ji} > 0}} ((hm_{ji} + d_{ji}^{(h)}) \log(hm_{ji} + d_{ji}^{(h)}) - (hm_{ji} + d_{ji}^{(h)}))$$

$$+ \sum_{(j,i) \in \mathcal{K}} (hm_{ji} + d_{ji}^{(h)}) \log \frac{\rho_i}{C_j} \Bigg]$$

$$= -\lim_{h \to \infty} \sum_{\substack{(j,i) \in \mathcal{K}: \\ m_{ji} > 0}} m_{ji} \log \frac{(m_{ji} + d_{ji}^{(h)}/h) C_j}{(m_j + d_j^{(h)}/h) \rho_i} = -\beta_\rho(m).$$

□

REMARK 3. The Kullback–Leibler divergence or relative entropy of distributions $p$ and $q$ on $\mathcal{I}$ is

$$D(p||q) = \sum_{i \in \mathcal{I}} p_i \log \left( \frac{p_i}{q_i} \right).$$

In our definition of $\beta_\rho(\cdot)$, if we define for each $j \in \mathcal{J}$, $p^j = (\frac{m_{ji}}{m_j} : i \ni j)$ and $q^j = (\frac{\rho_i}{\sum_{r \ni j} \rho_r} : i \ni j)$, then

$$\beta_\rho(m) = \sum_{j \,:\, m_j > 0} m_j D(p^j || q^j) + \sum_{j \,:\, m_j > 0} m_j \log\left(\frac{C_j}{\sum_{r \,:\, j \in r} \rho_r}\right).$$

So $\beta_\rho(\cdot)$ is a linear combination of Kullback–Leibler divergences. Normally we consider proportional fairness to maximize utility subject to constraints on flows. The duality given in Theorem 6.1 motivates us instead to view proportional fairness as minimizing entropy subject to constraints on packets.

Note that, since $x \log x$ is a continuous function, $\beta_\rho(\cdot)$ is a continuous function. We define $x \log x := 0$ when $x = 0$. Note also Lemma 6.1 applies for $d^{(h)} = hm - \lfloor hm \rfloor$. From this lemma the following result is reasonable.

PROPOSITION 6.1. *If $M$ is a random variable in $\mathbb{Z}_+^K$ with distribution (2.3), then, as $h \to \infty$, $\{\frac{M}{h}\}_{h \in \mathbb{N}}$ obeys a large deviation principle on $\mathbb{R}_+^K$ with good rate function $\beta_\rho(\cdot)$. That is for all $D \subset \mathbb{R}_+^K$ Borel measurable*

$$-\inf_{m \in D^\circ} \beta_\rho(m) \leq \liminf_{h \to \infty} \log \mathbb{P}\left(\frac{M}{h} \in D\right) \leq \limsup_{h \to \infty} \log \mathbb{P}\left(\frac{M}{h} \in D\right)$$
$$\leq -\inf_{m \in \bar{D}} \beta_\rho(m),$$

*where $D^\circ$ is the interior of $D$ and $\bar{D}$ is the closure of $D$.*

A proof of this proposition can be found in the Appendix. To prove the main theorem of this section we require two lemmas.

LEMMA 6.2. *For all $\Lambda \in (0, \infty)^I$,*

$$\inf_{m \in \mathbb{R}_+^K} \beta_\Lambda(m) = \begin{cases} 0, & \text{if } \sum_{i \,:\, j \in i} \Lambda_i \leq C_j, \forall j \in \mathcal{J}, \\ -\infty, & \text{otherwise.} \end{cases}$$

LEMMA 6.3. *For all $\Lambda \in (0, \infty)^I$, $\beta_\Lambda(\cdot)$ is a convex function.*

For proofs of Lemmas 6.2 and 6.3 see the Appendix. We now prove the main theorem of this section. We use the contraction principle (see [6], page 126). Recall that distribution (2.5) is the stationary distribution for the number of documents in transfer for all open multi-class queueing networks discussed in Section 2 and the spinning network. The following theorem expresses an important duality between network state and network flow.

THEOREM 6.1. *If $N$ is a random variable in $\mathbb{Z}_+^I$ with distribution (2.5) then as $h \to \infty$, $\{\frac{N}{h}\}_{h \in \mathbb{N}}$ obeys a large deviation principle on $\mathbb{R}_+^I$ with good rate function*

(6.3)
$$\alpha_\rho(n) := \min_{\substack{m \in \mathbb{R}_+^K}} \sum_{\substack{(j,i) \in \mathcal{K}: \\ m_j > 0}} m_{ji} \log \frac{m_{ji} C_j}{m_j \rho_i} \quad \text{subject to}$$
$$\sum_{j:\, j \in i} m_{ji} = n_i \quad \forall i \in \mathcal{I}$$

(6.4)
$$= \max_{\Lambda \in \mathbb{R}_+^I} \sum_{i \in \mathcal{I}} n_i \log \frac{\Lambda_i}{\rho_i} \quad \text{subject to} \quad \sum_{i:\, j \in i} \Lambda_i \leq C_j \quad \forall j \in \mathcal{J}.$$

*That is for all $A \subset \mathbb{R}_+^I$ Borel measurable, we have that*

$$-\inf_{n \in A^\circ} \alpha_\rho(n) \leq \liminf_{h \to \infty} \log \mathbb{P}\left(\frac{N}{h} \in A\right) \leq \limsup_{h \to \infty} \log \mathbb{P}\left(\frac{N}{h} \in A\right) \leq -\inf_{n \in \bar{A}} \alpha_\rho(n),$$

*where $A^\circ$ is the interior of $A$ and $\bar{A}$ is the closure of $A$.*

PROOF. Apply the contraction principle to Proposition 6.1 using continuous map $f : \mathbb{R}_+^K \to \mathbb{R}_+^I$ such that $f(m) = (\sum_{j:\, j \in i} m_{ji} : i \in \mathcal{I})$. This gives that $\{\frac{N}{h}\}_{h \in \mathbb{N}}$ obeys a large deviation principle with good rate function

$$\alpha_\rho(n) = \min_{\substack{m \in \mathbb{R}_+^K}} \sum_{\substack{(j,i) \in \mathcal{K}: \\ m_j > 0}} m_{ji} \log \frac{m_{ji} C_j}{m_j \rho_i} \quad \text{subject to} \quad \sum_{j:\, j \in i} m_{ji} = n_i \quad \forall i \in \mathcal{I}.$$

By Lemma 6.3, this is a convex optimization problem. Let us calculate its dual formulation. Taking Lagrange multipliers $\lambda \in \mathbb{R}^I$, its Lagrangian is

$$L(m, \lambda) = \sum_{\substack{(j,i) \in \mathcal{K}: \\ m_j > 0, n_i > 0}} m_{ji} \log \frac{m_{ji} C_j}{m_j \rho_i} + \sum_{i:\, n_i > 0} \lambda_i \left( n_i - \sum_{j:\, j \in i} m_{ji} \right)$$

$$= \sum_{\substack{(j,i) \in \mathcal{K}: \\ m_j > 0, n_i > 0}} m_{ji} \log \frac{m_{ji} C_j}{m_j \rho_i e^{\lambda_i}} + \sum_{i:\, n_i > 0} \lambda_i n_i.$$

By Lemma 6.2,

$$\min_{m \in \mathbb{R}_+^K} L(m, \lambda) = \begin{cases} \sum_{i:\, n_i > 0} n_i \lambda_i, & \text{if } \sum_{i:\, j \in i} \rho_i e^{\lambda_i} \leq C_j, \forall j \in \mathcal{J}, \\ -\infty, & \text{otherwise.} \end{cases}$$

Thus we find dual

$$\alpha_\rho(n) = \max_{\lambda \in \mathbb{R}^I} \sum_{i:\, n_i > 0} n_i \lambda_i \quad \text{subject to} \quad \sum_{i:\, j \in i} \rho_i e^{\lambda_i} \leq C_j.$$

Substituting $\Lambda_i = \rho_i e^{\lambda_i}$ gives

$$\alpha_\rho(n) = \max_{\Lambda \in \mathbb{R}_+^I} \sum_{i:\, n_i > 0} n_i \log \frac{\Lambda_i}{\rho_i} \quad \text{subject to} \quad \sum_{i:\, j \in i} \Lambda_i \leq C_j \quad \forall j \in \mathcal{J}. \qquad \square$$

**7. Teleological description.** In Section 6 we saw that a proportionally fair rate function determined the large deviations behavior of the stationary distribution of an open multi-class queueing network. In this section we discuss what this means for the behavior of packets and for the rate of document transfer in these networks. The queueing networks defined in Section 2 have no prescribed optimization structure. Even so, as the number of documents gets large, a network will be restricted to its most probable states, and because of this the network behaves as an optimizer. This notion of a queueing network collapsing to its most probable states is analogous to the heavy-traffic notion of state space collapse [4, 16, 24].

As in Section 6, we study the limit as the number of documents in transfer gets large but in proportion to some fixed vector. In this section, we characterize the solutions to the primal problem (6.3). In Theorem 7.1 show that the state of packets in an open multi-class queueing network with spinning converges in probability to the set of solutions of the primal problem. In Corollary 7.1, we show that at each queue the number of packets in transfer on each route converges in $L^1$ to a proportionally fair proportion of the number of packets at the queue. Finally in Theorem 7.2, we show that the stationary rate documents are transferred through these queueing networks converges to a proportionally fair allocation. We define

$$\beta(m) := \sum_{\substack{(j,i) \in \mathcal{K}: \\ m_j > 0}} m_{ji} \log \frac{m_{ji} C_j}{m_j} \quad \forall m \in \mathbb{R}_+^K.$$

We now characterize the solutions to the primal problem (6.3).

PROPOSITION 7.1. *Given $n \in \mathbb{R}_+^I$ and $m^* \in \mathbb{R}_+^K$ such that for all $i \in \mathcal{I}$, $\sum_{j \in i} m_{ji}^* = n_i$ then $m^*$ solves primal problem*

(7.1) $$\min_{m \in \mathbb{R}_+^K} \beta(m) \quad \text{subject to} \quad \sum_{j \in i} m_{ji} = n_i \quad \forall i \in \mathcal{I},$$

*if $\forall (j,i) \in \mathcal{K}$*

(7.2) $$m_{ji}^* C_j = m_j^* \Lambda_i^{\mathrm{PF}}(n).$$

PROOF. Suppose $m^*$ minimizes (7.1). We show (7.2) holds for some fixed $(j,i) \in \mathcal{K}$. If $m_j^* = 0$ then $m_{ji}^* = 0$ thus (7.2) holds. Now assume that $m_j^* > 0$, by the strong duality of primal (6.3) and dual (6.4) we have

$$\sum_{r:n_r>0} n_r \log \Lambda_r^{\mathrm{PF}}(n)$$

$$= \min_{m\in\mathbb{R}_+^K} \sum_{\substack{(l,r)\in\mathcal{K}:\\m_l>0,n_r>0}} m_{lr} \log \frac{m_{lr}}{m_l} C_l + \sum_{r:n_r>0} \log \Lambda_r^{\mathrm{PF}}(n)\left(n_r - \sum_{l:l\in r} m_{lr}\right)$$

and also $m^*$ minimizes this Lagrangian problem. The above expression gives

$$\min_{m\in\mathbb{R}_+^K} \sum_{l:m_l>0} m_l \sum_{\substack{r:l\in r,\\n_r>0}} \frac{m_{lr}}{m_l} \log \frac{m_{lr} C_l}{m_l \Lambda_r^{\mathrm{PF}}(n)} = 0.$$

By Lemma A.5, if $m_j^* > 0$ then $\forall i \ni j$, $m_{ji}^* C_j = m_j^* \Lambda_i^{\mathrm{PF}}(n)$, so (7.2) holds.

Now we prove the converse. If (7.2) holds then $\forall j \in \mathcal{J}$ such that $m_j^* > 0$

$$\frac{m_{ji}^* C_j}{m_j^*} = \Lambda_i^{\mathrm{PF}}(n),$$

therefore,

$$\beta(m^*) = \sum_{r\in\mathcal{I}} \sum_{l:l\in r} m_{lr}^* \log \Lambda_r^{\mathrm{PF}}(n) = \sum_{r:n_r>0} n_r \log \Lambda_r^{\mathrm{PF}}(n). \qquad \square$$

We define $\forall n \in \mathbb{R}_+^I$,

$$\mathcal{M}(n) = \left\{ m \in \mathbb{R}_+^K : m_{ji} C_j = m_j \Lambda_i^{\mathrm{PF}}(n) \ \forall (j,i) \in \mathcal{K}, \sum_{j:j\in i} m_{ji} = n_i \ \forall i \in \mathcal{I} \right\}.$$

In heavy-traffic literature $\mathcal{M}(n)$ would be thought of as an *invariant manifold*. Note that if a network of processor sharing queues were in state $m^* \in \mathcal{M}(n)$, then the transfer rate allocated to route $i$ packets would be $\Lambda_i^{\mathrm{PF}}(n)$. Since rare events occur in the most likely way, one may expect, as the number of documents in transfer gets large, that the state of the queues in the network will be close to $\mathcal{M}(n)$. Thus the rate of document transfer will be close to a proportionally fair allocation. We use this intuition to prove the next theorem but first we will require a lemma. Lemma 7.1 shows that the rate decay in probability is larger away from the manifold $\mathcal{M}(n)$.

LEMMA 7.1. *For all $\varepsilon > 0$, $\exists f(\varepsilon) > 0$ and $\delta(\varepsilon) > 0$ such that $\forall \delta < \delta(\varepsilon)$*

$$\sum_{i:n_i>0} n_i \log \frac{\Lambda_i^{\mathrm{PF}}(n)}{\rho_i} + f(\varepsilon) \leq \beta_{\varepsilon,\delta}^*.$$

*Where we define $\forall \varepsilon \geq 0$ and $\delta \geq 0$*

$$\beta_{\varepsilon,\delta}^* := \min_{m\in\mathbb{R}_+^K} \beta_\rho(m) \quad subject\ to \quad \max_{i\in\mathcal{I}} \left|\sum_{j\in i} m_{ji} - n_i\right| \leq \delta$$

*and*

$$\inf_{m' \in \mathcal{M}(n)} \|m - m'\| \geq \varepsilon.$$

See the Appendix for a proof of this result. We could interpret the following result as a state space collapse result [4, 16]; the result shows that the state of the queueing network converges in probability to the invariant manifold. Recall that all open queueing networks in Section 2 have stationary distribution (2.3).

THEOREM 7.1. *Let $M$ be any stochastic process on $\mathbb{R}_+^K$ with stationary distribution (2.3) and let $N_i = \sum_{j \in i} M_{ji}$, $\forall i \in \mathcal{I}$ then $\forall n \in \mathbb{R}_+^I$ and $\varepsilon > 0$*

$$\mathbb{P}\bigg(\inf_{m' \in \mathcal{M}(n)} \bigg\|\frac{M}{h} - m'\bigg\| \geq \varepsilon \bigg| N = \lfloor hn \rfloor\bigg) \xrightarrow[h \to \infty]{} 0.$$

PROOF. For all $n \in \mathbb{R}_+^I$ and $\forall \delta > 0$, one has that eventually in $h$,

$$\mathbb{P}\bigg(\inf_{m' \in \mathcal{M}(n)} \bigg\|\frac{M}{h} - m'\bigg\| \geq \varepsilon \bigg| N = \lfloor hn \rfloor\bigg)$$

$$= \frac{B}{B_{\lfloor hn \rfloor} \prod_{i \in \mathcal{I}} \rho_i^{\lfloor hn_i \rfloor}} \mathbb{P}\bigg(\inf_{m' \in \mathcal{M}(n)} \bigg\|\frac{M}{h} - m'\bigg\| \geq \varepsilon, N = \lfloor hn \rfloor\bigg)$$

(7.3)
$$\leq \frac{B}{B_{\lfloor hn \rfloor} \prod_{i \in \mathcal{I}} \rho_i^{\lfloor hn_i \rfloor}} \mathbb{P}\bigg(\inf_{m' \in \mathcal{M}(n)} \bigg\|\frac{M}{h} - m'\bigg\| \geq \varepsilon, \bigg\|\frac{N}{h} - n\bigg\| \leq \delta\bigg)$$

$$\leq \frac{B}{B_{\lfloor hn \rfloor} \prod_{i \in \mathcal{I}} \rho_i^{\lfloor hn_i \rfloor}} \exp\bigg\{-h\beta_{\varepsilon,\delta}^* + h\frac{f(\varepsilon)}{2}\bigg\}$$

$$\leq \frac{B}{B_{\lfloor hn \rfloor} \prod_{i \in \mathcal{I}} \rho_i^{\lfloor hn_i \rfloor}} \exp\bigg\{-h \sum_{i:n_i>0} n_i \log \frac{\Lambda_i^{\mathrm{PF}}(n)}{\rho_i} - h\frac{f(\varepsilon)}{2}\bigg\}.$$

We used Proposition 6.1 for the second inequality and Lemma 7.1 for the final inequality. Taking any $m^* \in \mathcal{M}(n)$, there exists a sequence $m^{(h)}$ with $hm^{(h)} \in \mathbb{Z}_+^K$, $\sum_{j:j \in i} hm_{ji}^{(h)} = \lfloor hn_i \rfloor$ $\forall i \in \mathcal{I}$ and $|hm_{ji}^{(h)} - hm_{ji}^*| \leq 2$ $\forall (j,i) \in \mathcal{K}$. So by Lemma 6.1 one has

$$\lim_{h \to \infty} \frac{1}{h} \log \mathbb{P}(M = hm^{(h)}) = -\beta_\rho(m^*) = -\sum_{i:n_i>0} n_i \log \frac{\Lambda_i^{\mathrm{PF}}(n)}{\rho_i}.$$

Thus $\exists h'$ such that $\forall h > h'$

(7.4)  $$\mathbb{P}(M = hm^{(h)}) \geq \exp\bigg\{-h \sum_{i:n_i>0} n_i \log \frac{\Lambda_i^{\mathrm{PF}}(n)}{\rho_i} - h\frac{f(\varepsilon)}{4}\bigg\}.$$

Hence one has that, eventually in $h$,

$$\mathbb{P}\left(\inf_{m' \in \mathcal{M}(n)} \left\| \frac{M}{h} - m' \right\| \geq \varepsilon \middle| N = \lfloor hn \rfloor \right)$$

$$\leq \frac{B}{B_{\lfloor hn \rfloor} \prod_{i \in \mathcal{I}} \rho_i^{\lfloor hn_i \rfloor}}$$

$$\times \exp\left\{ -h \sum_{i : n_i > 0} n_i \log \frac{\Lambda_i^{\mathrm{PF}}(n)}{\rho_i} - h \frac{f(\varepsilon)}{2} \right\}$$

$$\leq \frac{B}{B_{\lfloor hn \rfloor} \prod_{i \in \mathcal{I}} \rho_i^{\lfloor hn_i \rfloor}} \mathbb{P}(M = hm^{(h)}) e^{-hf(\varepsilon)/4}$$

$$= \mathbb{P}(M = hm^{(h)} | N = \lfloor hn \rfloor) e^{-hf(\varepsilon)/4} \xrightarrow[h \to \infty]{} 0.$$

We used (7.3) for the first inequality and (7.4) for the second inequality. □

REMARK 4. Pittel [20] proves Theorem 7.1 under the assumption that the manifold $\mathcal{M}(n)$ consists of a single point. This assumption meant Pittel did not require a result like Proposition 7.1. More generally Pittel's arguments do not make the primal–dual connection between queueing networks and proportional fairness. Other networks have a similar primal–dual large deviations connection between the network state and network flow (e.g., Loss Networks, see Kelly [12]).

We now consider what Theorem 7.1 implies for packets at each queue.

COROLLARY 7.1. $\forall n \in \mathbb{R}_+^I$ and $\forall (j, i) \in \mathcal{K}$

$$\frac{1}{h} \mathbb{E}[|M_{ji} C_j - M_j \Lambda_i^{\mathrm{PF}}(n)| | N = \lfloor hn \rfloor] \xrightarrow[h \to \infty]{} 0.$$

PROOF. Take $\varepsilon > 0$, by Theorem 7.1 and as $m_{ji} C_j = m_j \Lambda_i^{\mathrm{PF}}(n)$ $\forall m \in \mathcal{M}(n)$

$$\frac{1}{h} \mathbb{E}[|M_{ji} C_j - M_j \Lambda_i^{\mathrm{PF}}(n)| | N = \lfloor hn \rfloor]$$

$$= \frac{1}{h} \mathbb{E}\left[ |M_{ji} C_j - M_j \Lambda_i^{\mathrm{PF}}(n)| \mathbb{I}\left[ \inf_{m' \in \mathcal{M}(n)} \left\| \frac{M}{h} - m' \right\| \geq \varepsilon \right] \middle| N = \lfloor hn \rfloor \right]$$

$$+ \frac{1}{h} \mathbb{E}\left[ |M_{ji} C_j - M_j \Lambda_i^{\mathrm{PF}}(n)| \mathbb{I}\left[ \inf_{m' \in \mathcal{M}(n)} \left\| \frac{M}{h} - m' \right\| < \varepsilon \right] \middle| N = \lfloor hn \rfloor \right]$$

$$\leq 2 \left( \sum_{i : j \in i} n_i \right) C_j \mathbb{P}\left( \sup_{m' \in \mathcal{M}(n)} \left\| \frac{M}{h} - m' \right\| \geq \varepsilon \middle| N = \lfloor hn \rfloor \right) + (1 + I) C_j \varepsilon$$

$$\xrightarrow[h \to \infty]{} (1 + I) C_j \varepsilon.$$

For the first term before the inequality we use that $M_{ji} < M_j < \sum_{i \ni j} hn_i$ and $\Lambda_i^{\text{PF}}(n) < C_j$. For the second term we use that $|M_{ji} - hm'_{ji}| < h\varepsilon$ and that $m'$ satisfies (7.2). Since $\varepsilon$ is arbitrary the result holds. $\square$

The following result is the formal statement of (5.1) in Section 5.

LEMMA 7.2. *For all $(j,i) \in \mathcal{K}$ and for all $n \in \mathbb{Z}_+^I$ with $n_i > 0$,*

$$(7.5) \qquad \Lambda_i^{\text{SN}}(n)\mathbb{E}\left[\frac{M_j+1}{C_j}\bigg|N=n-e_i\right] = \mathbb{E}_n M_{ji}.$$

PROOF. Since $\mathbb{E}_{n-e_i}[\frac{M_j+1}{C_j}]$ is the expected sojourn of a route $i$ packet at queue $j$ the above result is really a statement of Little's law. We can show the result by explicit calculation,

$$\frac{\mathbb{E}_n M_{ji}}{\Lambda_i^{\text{SN}}(n)} = \frac{B_n}{B_{n-e_i}} \frac{1}{B_n} \sum_{\substack{m \in S(n): \\ m_{ji} > 0}} m_{ji} \prod_{l \in \mathcal{J}}\left(\binom{m_l}{m_{lr}:r \ni l}\prod_{r:l \in r}\left(\frac{1}{C_l}\right)^{m_{lr}}\right)$$

$$= \frac{1}{B_{n-e_i}} \sum_{\substack{m \in S(n): \\ m_{ji} > 0}} \frac{m_j}{C_j} \prod_{l \in \mathcal{J}}\left(\binom{m_l - \delta_{l,j}}{m_{lr} - \delta_{lr,ji}:r \ni l}\prod_{r:l \in r}\left(\frac{1}{C_l}\right)^{m_{lr}-\delta_{lr,ji}}\right)$$

$$= \frac{1}{B_{n-e_i}} \sum_{m' \in S(n-e_i)} \frac{m'_j+1}{C_j} \prod_{l \in \mathcal{J}}\left(\binom{m'_l}{m'_{lr}:r \ni l}\prod_{r:l \in r}\left(\frac{1}{C_l}\right)^{m'_{lr}}\right)$$

$$= \mathbb{E}_{n-e_i}\frac{M_j+1}{C_j}.$$

We define $\delta_{x,y} := 1$ for $x = y$ and $\delta_{x,y} := 0$ otherwise. Above we canceled terms and substituted $m'_{lr} = m_{lr} - \delta_{lr,ji}$. $\square$

Consider any of the multi-class queueing networks considered in Section 2. Recall from Corollary 2.3 that, given the number of documents in transfer is $n$, the stationary rate route $i$ documents are transferred through the network is $\Lambda_i^{\text{SN}}(n)$, the spinning allocation. We now prove that this rate converges to a proportionally fair bandwidth allocation. In this asymptotic sense these queueing networks behave as a proportionally fair optimizer.

THEOREM 7.2. *For all $n \in \mathbb{R}_+^I$ and $\forall i \in \mathcal{I}$*

$$\Lambda_i^{\text{SN}}(\lfloor hn \rfloor) \xrightarrow[h \to \infty]{} \Lambda_i^{\text{PF}}(n).$$

PROOF. If $n_i = 0$ the result is trivially true, so we assume $n_i > 0$. By the arguments in Theorem 7.1 and Corollary 7.1, one can see that

$$\frac{1}{h}\mathbb{E}[|M_{ji}C_j - M_j\Lambda_i^{\text{PF}}(n)||N = \lfloor hn \rfloor - e_i] \xrightarrow[h \to \infty]{} 0;$$

therefore,

$$\frac{1}{h}\mathbb{E}\left[M_{ji} - \frac{M_j}{C_j}\Lambda_i^{\text{PF}}(n)\Big|N = \lfloor hn \rfloor - e_i\right] \xrightarrow[h \to \infty]{} 0.$$

Summing over $j \in i$ gives

$$n_i - \frac{\Lambda_i^{\text{PF}}(n)}{h}\mathbb{E}\left[\sum_{j \in i}\frac{M_j}{C_j}\Big|N = \lfloor hn \rfloor - e_i\right] \xrightarrow[h \to \infty]{} 0;$$

substituting expression (7.5) gives

$$n_i - \frac{\Lambda_i^{\text{PF}}(n)}{h}\left(\frac{\lfloor n_i h \rfloor}{\Lambda_i^{\text{SN}}(\lfloor hn \rfloor)} - \sum_{j \in i}\frac{1}{C_j}\right) \xrightarrow[h \to \infty]{} 0$$

and finally rearranging gives the result

$$\Lambda_i^{\text{SN}}(\lfloor hn \rfloor) \xrightarrow[h \to \infty]{} \Lambda_i^{\text{PF}}(n). \qquad \square$$

## APPENDIX

**A.1. Proof of convergence to the spinning network.** The proof provided here gives the first rigorous proof of a packet level model of document transfer converging weakly to a stochastic flow level model.

In this section and as defined in Section 3, for $c \in \mathbb{N}$, let $M^{(c)}$ be the number of packets in transfer on each route and at each queue in an open multi-class queueing network with spinning; let $N^{(c)}$ be the number of documents in transfer in this queueing network and let $N^{(\infty)}$ be the spinning network. Finally for $c \in \mathbb{N} \cup \{\infty\}$, let $\tau^{k,(c)}$ be the $k$th jump of $N^{(c)}$, that is, the time of the $k$th document arrival or departure event. We set $\tau^{0,(c)} := 0$.

To prove Theorem 3.1 we use a coupling argument. We prove under this coupling, as $c \to \infty$, $N^{(c)}$ eventually makes the same jumps as $N^{(\infty)}$ and that the associated jump times converge almost surely. This is sufficient to prove convergence in the Skorohod topology:

LEMMA A.1. *For $c \in \mathbb{N} \cup \{\infty\}$, let $N^{(c)}:[0, \infty) \to \mathbb{Z}_+^I$ be a nonexplosive jump processes with increasing jump times $(\tau^{k,(c)} : k \in \mathbb{Z}_+)$. If for all $k \in \mathbb{Z}_+$*

$$\tau^{k,(c)} \xrightarrow[c \to \infty]{} \tau^{k,(\infty)} \quad \text{and} \quad N^{(c)}(\tau^{k,(c)}) \xrightarrow[c \to \infty]{} N^{(\infty)}(\tau^{k,(\infty)}),$$

*then in the Skorohod topology on $[0, 1]$,*

$$N^{(c)} \xrightarrow[c \to \infty]{} N^{(\infty)}.$$

For a proof of this result see Billingsley [1], page 137. We now work to form a coupling so that the first jump and jump time converge. That is, we will prove:

PROPOSITION A.1. *Let $n^0 \in \mathbb{Z}_+^I$ and for each $c \in \mathbb{N}$ take a state $m^{0,(c)} \in S(n^0)$. There exists a coupling of $N^{(\infty)}$ with initial position $n_0$ and $M^{(c)}$ with initial position $m^{0,(c)}$ such that, almost surely*

$$\tau^{1,(c)} \xrightarrow[c \to \infty]{} \tau^{1,(\infty)} \quad and \quad N^{(c)}(\tau^{1,(c)}) \xrightarrow[c \to \infty]{} N^{(\infty)}(\tau^{1,(\infty)}).$$

To prove Theorem 3.1 we will apply this result to each jump interval of $N^{(c)}$. The proof of Proposition A.1 couples each $M^{(c)}$ with a single closed queueing network. The result is then an application of the renewal theorem.

As described in Section 2, let $\bar{M}$ be the number of packets in transfer on each route and at each queue for a closed multi-class queueing network with $n^0 \in \mathbb{Z}_+^I$ documents in transfer and service capacities $(C_j : j \in \mathcal{J})$. For states $m \in S(n^0)$, define $\sigma_m$ to be the first time $\bar{M}$ visits state $m$. Since $\bar{M}$ is recurrent, almost surely, $\sigma_m < \infty$. As noted in Section 3, $M^{(c)}$ will behave as a closed queueing network until the first document arrival or departure time. In particular, we will define

$$(A.1) \qquad M^{(c)}(t) := \bar{M}(ct + \sigma_{m^{0,(c)}}) \qquad \forall t \in [0, \tau^{1,(c)}).$$

The $ct$ term ensures the correct transition rates and the $\sigma_{m^{0,(c)}}$ ensures $M^{(c)}$ has the correct initial state. We will formally define $\tau^{1,(c)}$ later.

Let $D_i(t)$ be the number of route $i$ packets to have been served at the final queue of route $i$ in closed queueing network $\bar{M}$ by time $t$. By Corollary 2.3, we know the stationary throughput of route $i$ packets at any queue $j \in i$ is $\Lambda_i^{\mathrm{SN}}(n^0)$. Thus we can prove the following renewal lemma and corollary.

LEMMA A.2. *Almost surely, for all $\eta > 0$,*

$$\sup_{t \in [0,\eta]} \left| \frac{D_i(ct)}{c} - \Lambda_i^{\mathrm{SN}}(n^0) t \right| \xrightarrow[c \to \infty]{} 0.$$

PROOF. Let $R_{i,m}(t)$ be the number of route $i$ packets to have been served at the final queue on route $i$ by the closed queueing network $\bar{M}$ when it is in state $m$. Let $\gamma_i(m)$ be the drift of $R_{i,m}$. For any Markov chain the process that records the current state of the Markov chain and the next state is also a Markov chain. So $R_{i,m}$ is a renewal process and thus obeys the functional renewal theorem. That is, almost surely, for all $\eta > 0$,

$$\sup_{t \in [0,\eta]} \left| \frac{R_{i,m}(ct)}{c} - \gamma_i(m) t \right| \xrightarrow[c \to \infty]{} 0.$$

For a proof of this see Chen and Yao [5], page 106. By the definition of $D_i(t)$ and Corollary 2.3, we know that

$$D_i(t) = \sum_{m \in S(n^0)} R_{i,m}(t) \quad \text{and} \quad \Lambda_i^{\text{SN}}(n^0) = \sum_{m \in S(n^0)} \gamma_i(m).$$

So, almost surely, $\forall \eta > 0$

$$\sup_{t \in [0,\eta]} \left| \frac{D_i(ct)}{c} - \Lambda_i^{\text{SN}}(n^0)t \right| \leq \sum_{m \in S(n^0)} \sup_{t \in [0,\eta]} \left| \frac{R_{i,m}(ct)}{c} - \gamma_i(m)t \right| \xrightarrow[c \to \infty]{} 0. \quad \square$$

COROLLARY A.2. *Almost surely,*

$$\sup_{t \in [0,\eta]} \left| \frac{D_i(ct + \sigma_{m^{0,(c)}})}{c} - \Lambda_i^{\text{SN}}(n^0)t \right| \xrightarrow[c \to \infty]{} 0.$$

PROOF. As $\bar{M}$ is recurrent on all states in $S(n^0)$, almost surely, $\sigma_m < \infty$ $\forall m \in S(n^0)$. Thus, by this and Lemma A.2, almost surely,

$$\sup_{t \in [0,\eta]} \left| \frac{D_i(ct + \sigma_{m^{0,(c)}})}{c} - \Lambda_i^{\text{SN}}(n^0)\left(t + \frac{\sigma_{m^{0,(c)}}}{c} - \frac{\sigma_{m^{0,(c)}}}{c}\right) \right|$$

$$< \Lambda_i^{\text{SN}}(n^0) \frac{\sigma_{m^{0,(c)}}}{c} + \sup_{t \in [0,\eta+\sigma_{m^{0,(c)}}]} \left| \frac{D_i(ct)}{c} - \Lambda_i^{\text{SN}}(n^0)t \right| \xrightarrow[c \to \infty]{} 0. \quad \square$$

In the open queueing network $M^{(c)}$, we suppose each document in transfer on route $i$ is geometrically distributed with parameter $\frac{\mu_i}{c}$. Thus, if $n_i^0 > 0$ and assuming no other document arrival or departures occur, the time until the first route $i$ document is transferred is

$$S_i^{(c)} := \inf\{t : D_i(ct + \sigma_{m^{0,(c)}}) = Y_i^{(c)}\},$$

where $Y_i^{(c)}$ is geometrically distributed with parameter $\frac{\mu_i}{c}$ and $(Y_i^{(c)} : c \in \mathbb{N})$ is independent of $\bar{M}$. Now suppose

(A.2) $$\frac{Y_i^{(c)}}{c} \xrightarrow[c \to \infty]{} Y_i^{(\infty)} \quad \text{almost surely,}$$

where $Y_i^{(\infty)}$ is exponentially distributed with parameter $\mu_i$. This can be constructed using the Skorohod representation theorem. Using this and Corollary A.2 we can show that $S_i^{(c)}$ converges to an exponential distribution.

LEMMA A.3. *For each $i \in \mathcal{I}$ such that $n_i > 0$, almost surely*

$$S_i^{(c)} \xrightarrow[c \to \infty]{} S_i^{(\infty)},$$

*where $S_i^{(\infty)}$ is exponentially distributed with parameter $\mu_i \Lambda_i^{\text{SN}}(n^0)$.*

PROOF. Define $S_i^{(\infty)} = \frac{Y_i^{(\infty)}}{\Lambda_i^{\text{SN}}(n^0)}$. By Corollary A.2 and (A.2), almost surely, $\forall \varepsilon > 0$ and $\forall \eta > \frac{Y_i^{(\infty)} + 2\varepsilon}{\Lambda_i^{\text{SN}}(n^0)}$, $\exists c'$ such that $\forall c > c'$,

$$\sup_{t \in [0,\eta]} \left| \frac{D_i(ct + \sigma_{m^{0,(c)}})}{c} - \Lambda_i^{\text{SN}}(n^0) t \right| < \varepsilon, \qquad \left| \frac{Y_i^{(c)}}{c} - Y_i^{(\infty)} \right| < \varepsilon.$$

Hence

$$\frac{1}{c} D_i \left( \frac{cY_i^{(\infty)}}{\Lambda_i^{\text{SN}}(n^0)} - \frac{2c\varepsilon}{\Lambda_i^{\text{SN}}(n^0)} + \sigma_{m^{0,(c)}} \right) \leq Y_i^{(\infty)} - \varepsilon < \frac{Y_i^{(c)}}{c},$$

thus

$$S_i^{(c)} = \inf\{t \geq 0 : D_i(ct + \sigma_{m^{0,(c)}}) = Y^{(c)}\} > \frac{Y_i^{(\infty)}}{\Lambda_i^{\text{SN}}(n^0)} - \frac{2\varepsilon}{\Lambda_i^{\text{SN}}(n^0)}$$

$$= S_i^{(\infty)} - \frac{2\varepsilon}{\Lambda_i^{\text{SN}}(n^0)}.$$

By a similar argument one can see that

$$S_i^{(c)} < S_i^{(\infty)} + \frac{2\varepsilon}{\Lambda_i^{\text{SN}}(n^0)}.$$

Thus $S_i^{(c)} \to S_i^{(\infty)}$ as $c \to \infty$, almost surely. □

We may choose $(Y_i^{(c)} : c \in \mathbb{N} \cup \{\infty\})$ independently for each $i \in \mathcal{I}$ with $n_i > 0$. The transfer of a document on route $i$ could be interrupted by an earlier document arrival or departure event. Thus letting $E_i$ be independent exponentially distributed with parameter $\nu_i$ for $i \in \mathcal{I}$. We are interested in the time when the first arrival or departure time occurs, so we consider

$$\tau^{1,(c)} := \min\{S_i^{(c)} : n_i > 0\} \wedge \min\{E_i : i \in \mathcal{I}\} \qquad \forall c \in \mathbb{N} \cup \{\infty\}.$$

By the last lemma we know that $\tau^{1,(c)} \to \tau^{1,(\infty)}$, as $c \to \infty$. The term achieving these minima determines which arrival or departure occurs. So we may define our coupled process, $M^{(c)}$, up until the first arrival or departure time by (A.1) and by

$$M^{(c)}(\tau^{1,(c)}) := \begin{cases} M^{(c)}(\tau^{1,(c)}-) + e_{j_0^i i}, & \text{if } \tau^{1,(c)} = E_i, \\ M^{(c)}(\tau^{1,(c)}-) - e_{j(i)i}, & \text{if } \tau^{1,(c)} = S_i^{(c)}. \end{cases}$$

The term $e_{ji} \in \mathbb{R}_+^K$ is the unit vector in the $(j,i) \in \mathcal{K}$ direction, thus $e_{j_0^i i}$ corresponds to the arrival of the first packet in a route $i$ document and $-e_{j(i)i}$ corresponds to the departure of the final packet in a route $i$ document. Also we may define $N^{(\infty)}$ by $N^{(\infty)}(t) := n^0$ for $t < \tau^{1,(\infty)}$ and by

$$N^{(\infty)}(\tau^{1,(\infty)}) := \begin{cases} n^0 + e_i, & \text{if } \tau^{1,(\infty)} = E_i, \\ n^0 - e_i, & \text{if } n_i > 0 \text{ and } \tau^{1,(\infty)} = S_i^{(\infty)}. \end{cases}$$

Similarly, $e_i \in \mathbb{R}_+^I$ is the $i$th unit vector. This defines our coupled process up to the first document arrival or departure time. We can now prove Proposition A.1.

PROOF OF PROPOSITION A.1.   By Lemma A.3 we know that, almost surely,

$$(\text{A.3}) \qquad \tau^{1,(c)} \xrightarrow[c \to \infty]{} \tau^{1,(\infty)}.$$

Also as $\{S_i^{(\infty)} : n_i > 0\} \cup \{E_i : i \in \mathcal{I}\}$ are independent exponentially distributed random variables, almost surely, no two terms are equal. Thus due to this and (A.3), almost surely, eventually as $c \to \infty$,

$$\arg\min[\{S_i^{(c)} : n_i > 0\} \cup \{E_i : i \in \mathcal{I}\}]$$
$$= \arg\min[\{S_i^{(\infty)} : n_i > 0\} \cup \{E_i : i \in \mathcal{I}\}],$$

and so, almost surely, as $c \to \infty$, $N^{(c)}(\tau^{1,(c)}) \to N^{(\infty)}(\tau^{1,(\infty)})$. □

Proposition A.1 guarantees convergence up to and including the first document arrival or departure time. Essentially repeating this argument constructs each interval $[\tau^{k,(c)}, \tau^{k+1,(c)}]$ and proves Theorem 3.1.

PROOF OF THEOREM 3.1.   On a single probability space we will inductively construct $N^{(\infty)}$ and $M^{(c)}$ for $c \in \mathbb{N}$. On the probability space we will prove that the following induction hypothesis holds $\forall \kappa \in \mathbb{Z}_+$: there exists a coupling of $M^{(c)}$ for $c \in \mathbb{N}$ and $N^{(\infty)}$ up to and including the $\kappa$th document's arrival or departure time such that, almost surely, $\forall k \leq \kappa$,

$$(\text{A.4}) \qquad \tau^{k,(c)} \xrightarrow[c \to \infty]{} \tau^{k,(\infty)} \quad \text{and} \quad N^{(c)}(\tau^{k,(c)}) \xrightarrow[c \to \infty]{} N^{(\infty)}(\tau^{k,(\infty)}).$$

Let us prove that the induction hypothesis holds for the case of $\kappa = 0$. We must find a coupling of $M^{(c)}(0)$ and $N^{(\infty)}(0)$ so that (A.4) holds. We know $\tau^{1,(c)} \to \tau^{1,(\infty)}$ as $c \to \infty$ holds as $\tau^{0,(c)} := 0$. By assumption $N^{(c)}(0)$ converges weakly to $N^{(\infty)}(0)$, so by the Skorohod representation theorem, we may choose $(N^{(c)}(0) : c \in \mathbb{N} \cup \{\infty\})$, such that, almost surely,

$$N^{(c)}(0) \xrightarrow[c \to \infty]{} N^{(\infty)}(0).$$

Given $N^{(c)}$ we know the required distribution of $M^{(c)}(0)$, so take $f^{(c)} : \mathbb{Z}_+^I \times [0,1] \to \mathbb{Z}_+^K$ such that for a uniform random variable $U$,

$$\mathbb{P}(f^{(c)}(n,U) = m) = \mathbb{P}(M^{(c)} = m | N^{(c)} = n) \qquad \forall m \in \mathbb{Z}_+^K, n \in \mathbb{Z}_+^I.$$

Therefore taking an independent uniform random variable we may define

$$M^{(c)}(0) = f^{(c)}(N^{(c)}(0), U) \qquad \forall c \in \mathbb{N}.$$

Thus $M^{(c)}(0)$ is of the correct distribution and (A.4) holds. This proves the induction hypothesis is true for $\kappa = 0$.

Suppose the induction hypothesis holds for $\kappa - 1$. Using Proposition A.1, we will show the induction hypothesis holds for $\kappa$ by extending the process $M^{(c)}$ from time $\tau^{\kappa-1,(c)}$ to the next document arrival or departure time $\tau^{\kappa,(c)}$.

By the induction hypothesis $\exists c' \in \mathbb{N}$ such that $\forall c > c'$, $N^{(c)}(\tau^{\kappa-1,(c)}) = N^{(\infty)}(\tau^{\kappa-1,(\infty)})$. Define

$$m^{\kappa-1,(c)} = M^{(c)}(\tau^{\kappa-1,(c)}) \qquad \forall c > c'.$$

Now apply Proposition A.1 with initial state $(m^{\kappa-1,(c)} : c > c')$ to give process $M^{\kappa,(c)}$ defined until the first document arrival or departure time $\tau^{k,(c)} - \tau^{k-1,(c)}$. We extend $M^{(c)}$ to include $t \in [\tau^{\kappa-1,(c)}, \tau^{\kappa,(c)}]$ by defining

$$M^{(c)}(t) := M^{\kappa,(c)}(t - \tau^{\kappa,(c)}) \qquad \text{for } t \in [\tau^{\kappa-1,(c)}, \tau^{\kappa,(c)}].$$

By Proposition A.1 we know (A.4) holds for $\kappa$. This completes the induction step. Thus the induction hypothesis holds for all $\kappa \in \mathbb{N}$. At each induction step we required a countable collection of independent random variables, thus the coupled processes $M^{(c)}$ and $N^{(\infty)}$ can be constructed on a probability space consisting of a countable set of independent random variables. Since (A.4) holds for all $k \in \mathbb{N}$, Lemma A.1 gives that, almost surely,

$$N^{(c)} \xrightarrow[c \to \infty]{} N^{(\infty)}$$

in the Skorohod topology on $[0,1]$. Thus this implies weak convergence. □

**A.2. Proof of Proposition 6.1.** To prove Proposition 6.1 we expand on Lemma 6.1 to make a large deviation principle. The following lemmas will be used. The proofs of both lemmas are simple calculus arguments.

LEMMA A.4. *For $m \in \mathbb{R}_+^K$ and $j \in \mathcal{J}$ fixed with $m_j > 0$. Maximizing over vectors $\theta = (\theta_i : i \ni j) > 0$ satisfying constraint $\sum_{i : j \in i} \rho_i e^{\theta_i} = C_j$, we have that*

$$\max_{\theta > 0} \left\{ \sum_{i : j \in i} \theta_i m_{ji} : \sum_{i : j \in i} \rho_i e^{\theta_i} = C_j \right\} = \sum_{i : j \in i} m_{ji} \log \frac{m_{ji} C_j}{m_j \rho_i}.$$

LEMMA A.5. *For fixed parameters $\Lambda \in (0,\infty)^I$ and $j \in \mathcal{J}$. Maximizing over probability distributions $p = (p_i : i \ni j)$, we have that*

$$\min_{p>0}\left\{\sum_{i\,:\,j\in i} p_i \log \frac{p_i C_j}{\Lambda_i} : \sum_{i\,:\,j\in i} p_i = 1\right\} = \log \frac{C_j}{\sum_{i\,:\,j\in i} \Lambda_i}$$

*and the minimum is attained by $p_i = \frac{\Lambda_i}{\sum_{r\,:\,j\in r} \Lambda_r}$.*

We now prove Proposition 6.1. The result makes use of the upper bound of the Gärtner–Ellis theorem (see [6], page 44).

PROOF OF PROPOSITION 6.1. First we prove the lower bound. Take $O \subset \mathbb{R}_+^K$ open, for all $m \in O$ $\exists h'$ such that $\forall h > h'$, $\frac{\lfloor hm \rfloor}{h} \in O$. By this and Lemma 6.1 we have that

$$\liminf_{h\to\infty} \frac{1}{h} \log \mathbb{P}\left(\frac{M}{h} \in O\right) \geq \sup_{m\in O} \liminf_{h\to\infty} \frac{1}{h} \log \mathbb{P}(M = \lfloor hm \rfloor) = -\inf_{m\in O} \beta_\rho(m).$$

This gives the required lower bound.

Now we prove the upper bound. $M$ has moment generating function

$$\mathbb{E} e^{\sum_{(j,i)\in\mathcal{K}} \theta_{ji} M_{ji}}$$
$$= \begin{cases} \prod_{j\in\mathcal{J}} \left(\dfrac{C_j - \sum_{i\,:\,j\in i} \rho_i}{C_j - \sum_{i\,:\,j\in i} \rho_i e^{\theta_{ji}}}\right), & \text{if } \sum_{i\,:\,j\in i} \rho_i e^{\theta_{ji}} < C_j, \forall j \in \mathcal{J}, \\ \infty, & \text{otherwise.} \end{cases}$$

Thus

$$F(\theta) := \lim_{h\to\infty} \frac{1}{h} \log \mathbb{E} e^{\sum_{(j,i)\in\mathcal{K}} h\theta_{ji} M_{ji}/h}$$
$$= \begin{cases} 0, & \text{if } \sum_{i\,:\,j\in i} \rho_i e^{\theta_{ji}} < C_j, \forall j \in \mathcal{J}, \\ \infty, & \text{otherwise.} \end{cases}$$

Let $F^*(\cdot)$ be the Legendre–Fenchel transform of $F(\cdot)$. By Lemma A.4,

$$F^*(m) = \sup_{\theta\in\mathbb{R}^K} \left\{\sum_{(j,i)\in\mathcal{K}} \theta_{ji} m_{ji} : \sum_{i\,:\,j\in i} \rho_i e^{\theta_{ji}} < C_j, \forall j \in \mathcal{J}\right\}$$
$$= \sum_{j\in\mathcal{J}} \sup_{\substack{\phi_i>0:\\i\ni j}} \left\{\sum_{i\,:\,j\in i} \phi_i m_{ji} : \sum_{i\,:\,j\in i} \rho_i e^{\phi_i} < C_j\right\}$$
$$= \sum_{j\,:\,m_j>0} \sum_{i\,:\,j\in i} m_{ji} \log \frac{m_{ji} C_j}{m_j \rho_i} = \beta_\rho(m),$$

$\forall m \in \mathbb{R}_+^K$. Thus the Gärtner–Ellis theorem gives for all closed sets $C \subset \mathbb{R}_+^K$

$$\limsup_{h\to\infty} \frac{1}{h} \log \mathbb{P}\left(\frac{M}{h} \in C\right) \leq -\inf_{m \in C} \beta_\rho(m).$$

Finally, we prove $\beta_\rho(\cdot)$ is a good rate function. $\beta_\rho(\cdot)$ is continuous with values in $\mathbb{R}_+$ and so is a rate function. By Lemma A.5, $\forall \alpha \geq 0$, $\forall m \in \mathbb{R}_+^K$ if

$$m_j > \frac{\alpha}{\log C_j - \log \sum_{i\,:\,j\in i} \rho_i},$$

then

$$\beta_\rho(m) = \sum_{j\,:\,m_j > 0} m_j \sum_{i\,:\,j\in i} \frac{m_{ji}}{m_j} \log \frac{m_{ji} C_j}{m_j \rho_i} \geq \sum_{j\,:\,m_j > 0} m_j \log \frac{C_j}{\sum_{i\,:\,j\in i} \rho_i} > \alpha.$$

Thus

$$\{m \in \mathbb{R}_+^K : \beta_\rho(m) \leq \alpha\} \subset \left\{m \in \mathbb{R}_+^K : 0 \leq m_j \leq \frac{\alpha}{\log C_j - \log \sum_{i\,:\,j\in i} \rho_i}\right\}.$$

So all level sets are compact and hence $\beta_\rho(\cdot)$ is a good rate function. $\square$

### A.3. Proofs of additional lemmas.

PROOF OF LEMMA 6.2. By Lemma A.5

$$\inf_{m \in \mathbb{R}_+^K} \beta_\Lambda(m) = \inf_{m \in \mathbb{R}_+^K} \sum_{j\,:\,m_j > 0} m_j \sum_{i\,:\,j\in i} \frac{m_{ji}}{m_j} \log \frac{m_{ji} C_j}{m_j \Lambda_i}$$

$$= \inf_{m' \in \mathbb{R}_+^J} \sum_{j\,:\,m'_j > 0} m'_j \log \frac{C_j}{\sum_{i\,:\,j\in i} \Lambda_i}$$

$$= \begin{cases} 0, & \text{if } \sum_{i\,:\,j\in i} \Lambda_i \leq C_j, \\ -\infty, & \text{otherwise.} \end{cases} \qquad \square$$

PROOF OF LEMMA 6.3. We first show $\beta_\Lambda(\cdot)$ is convex on $(0,\infty)^K$. One can see $\forall m \in (0,\infty)^K$ and $\forall (j,i), (l,r) \in \mathcal{K}$ that

$$\frac{\partial^2 \beta_\Lambda}{\partial m_{ji}\, \partial m_{lr}}(m) = \frac{1}{m_{ji}} \mathbb{I}[(j,i) = (l,r)] - \frac{1}{m_j} \mathbb{I}[l = r].$$

So for all vectors $a \in \mathbb{R}^K \setminus \{0\}$ the second partial derivative in direction $a$ is

$$\frac{\partial^2 \beta_\Lambda}{\partial a^2} = \sum_{\substack{(j,i)\in\mathcal{K},\\(l,r)\in\mathcal{K}}} a_{ji} a_{lr} \frac{\partial^2 \beta_\Lambda}{\partial m_{ji}\, \partial m_{lr}}(m) = \sum_{(j,i)\in\mathcal{K}} \frac{a_{ji}^2}{m_{ji}} - \sum_{j\in\mathcal{J}} \sum_{i\,:\,j\in i} \sum_{r\,:\,j\in r} \frac{a_{ji} a_{jr}}{m_j}.$$

We will prove that this expression is positive. Note that for fixed $j \in \mathcal{J}$,

$$\sum_{i:\, j \in i} a_{ji}^2 \frac{m_j}{m_{ji}} = \sum_{i:\, j \in i} \sum_{r:\, j \in r} \frac{1}{2}\left(a_{ji}^2 \frac{m_{jr}}{m_{ji}} + a_{jr}^2 \frac{m_{ji}}{m_{jr}}\right).$$

For all $i \neq r$, $a_{ji}^2 x + a_{jr}^2 \frac{1}{x}$ is convex on $(0, \infty)$ and has its minimum at $x = |\frac{a_{jr}}{a_{ji}}|$. So

$$\sum_{i:\, j \in i} a_{ji}^2 \frac{m_j}{m_{ji}} \geq \sum_{i:\, j \in i} \sum_{r:\, j \in r} |a_{ji}||a_{jr}| \geq \sum_{i:\, j \in i} \sum_{r:\, j \in r} a_{ji} a_{jr} \geq 0.$$

Dividing by $m_j$ and summing over $j \in \mathcal{J}$ gives that

$$\frac{\partial^2 \beta_\Lambda}{\partial a^2} = \sum_{j \in \mathcal{J}} \left(\sum_{i:\, j \in i} \frac{a_{ji}^2}{m_{ji}} - \sum_{i:\, j \in i} \sum_{r:\, j \in r} \frac{a_{ji} a_{jr}}{m_j}\right) \geq 0.$$

This proves $\beta_\Lambda(\cdot)$ is convex on $(0, \infty)^K$.

Take $m, \bar{m} \in \mathbb{R}_+^K$ and two sequences, $m^{(h)}, \bar{m}^{(h)} \in (0, \infty)^K$, that converge to $m$ and $\bar{m}$, respectively. By continuity of $\beta_\Lambda(\cdot)$, $\forall \theta \in (0, 1)$

$$\begin{aligned}\beta_\Lambda(\theta m + (1-\theta)\bar{m}) &= \lim_{h \to \infty} \beta_\Lambda(\theta m^{(h)} + (1-\theta)\bar{m}^{(h)}) \\ &\leq \lim_{h \to \infty} \theta \beta_\Lambda(m^{(h)}) + (1-\theta)\beta_\Lambda(\bar{m}^{(h)}) \\ &= \theta \beta_\Lambda(m) + (1-\theta)\beta_\Lambda(\bar{m}),\end{aligned}$$

thus $\beta_\Lambda(\cdot)$ is convex on $\mathbb{R}_+^K$. $\square$

PROOF OF LEMMA 7.1. It is clear that $\beta_{\varepsilon,\delta}^*$ is nonincreasing in $\delta$. We now claim that $\beta_{\varepsilon,\delta}^* \nearrow \beta_{\varepsilon,0}^*$ as $\delta \searrow 0$. If this were not so by compactness of our optimization region, we could choose $\delta_k \searrow 0$ and $m^k \in \mathbb{R}_+^K$, with

$$\max_{i \in \mathcal{I}}\left|\sum_{j \in i} m_{ji}^k - n_i\right| \leq \delta_k, \qquad \inf_{m' \in \mathcal{M}(n)} \|m^k - m'\| \geq \varepsilon, \qquad \sup_k \beta_\rho(m_k) < \beta_{\varepsilon,0}^*,$$

and such that for some $m \in \mathbb{R}_+^K$, $m^k \to m$ as $k \to \infty$. But then by the continuity of $\beta_\rho$, $\beta_\rho(m) < \beta_{\varepsilon,0}^*$. This then contradicts the minimality of $\beta_{\varepsilon,0}^*$. So, $\beta_{\varepsilon,\delta}^* \nearrow \beta_{\varepsilon,0}^*$ as $\delta \searrow 0$.

By strong duality of (6.3) and (6.4), $\beta_{0,0}^* = \sum_{i:\, n_i > 0} n_i \log \frac{\Lambda_i^{\mathrm{PF}}(n)}{\rho_i}$. As $\beta_\rho$ is continuous and $\mathcal{M}(n)$ is compact, $\beta_{\varepsilon,0}^* > \beta_{0,0}^*$ $\forall \varepsilon > 0$. Take any $f(\varepsilon) < \beta_{\varepsilon,0}^* - \beta_{0,0}^*$. By the last paragraph, $\exists \delta(\varepsilon)$ such that $\forall \delta < \delta(\varepsilon)$, $\beta_{\varepsilon,0}^* - \beta_{\varepsilon,\delta}^* < 2f(\varepsilon)$. This gives that $\beta_{\varepsilon,0}^* - \beta_{\varepsilon,\delta}^* < f(\varepsilon) + \beta_{\varepsilon,0}^* - \beta_{0,0}^*$, and hence the result. $\square$

**Acknowledgments.** I am most grateful to Frank Kelly whose insight and enthusiasm has appreciably aided this paper. I would like to thank Alexandre Proutière for various comments that have improved the content and context of this paper.

CENTRE FOR MATHEMATICAL SCIENCES
UNIVERSITY OF CAMBRIDGE
WILBERFORCE ROAD
CAMBRIDGE CB3 0WB
UNITED KINGDOM
E-MAIL: n.s.walton@statslab.cam.ac.uk
URL: http://www.statslab.cam.ac.uk/~nsw26